\documentclass[conference]{IEEEtran}

\usepackage{ifpdf}

\usepackage{url}
\ifpdf
	\usepackage[pdftex]{graphicx}
\else
	\usepackage[dvips]{graphicx}
\fi

\usepackage{url}
\usepackage[noadjust]{cite}
\DeclareGraphicsExtensions{.eps}
\usepackage{epstopdf}
\usepackage[font=footnotesize]{subfig}
\usepackage{caption}
\usepackage{psfrag}
\usepackage{amsmath}
\usepackage{booktabs,caption, dblfloatfix}

\newcommand\Tstrut{\rule{0pt}{2.6ex}}         % = `top' strut
\newcommand\Bstrut{\rule[-0.9ex]{0pt}{0pt}}   % = `bottom' strut

\begin{document}

\title{A Comparison of Policies on the Participation of Storage in U.S. Frequency Regulation Markets }

\author{
	\IEEEauthorblockN{Bolun Xu, \emph{Student Member, IEEE}, Yury Dvorkin, \emph{Student Member, IEEE}, Daniel S. Kirschen, \emph{Fellow, IEEE}}
	\IEEEauthorblockN{C. A. Silva-Monroy, \emph{Member, IEEE}, Jean-Paul Watson, \emph{Member, IEEE}}
		\IEEEauthorblockA{Emails: \{xubolun, dvorkin, kirschen\}@uw.edu, \{casilv, jwatson\}@sandia.gov}
}

\maketitle

\begin{abstract}

Because energy storage systems have better ramping characteristics than traditional generators, their participation in frequency regulation should facilitate the balancing of load and generation. However, they cannot sustain their output indefinitely. System operators have therefore implemented new frequency regulation policies to take advantage of the fast ramps that energy storage systems can deliver while alleviating the problems associated with their limited energy capacity. This paper contrasts several U.S. policies that directly affect the participation of energy storage systems in frequency regulation and compares the revenues that the owners of such systems might achieve under each policy.
\end{abstract}


\IEEEpeerreviewmaketitle

\section{Introduction}
\let\thefootnote\relax\footnote{This paper will be presented in IEEE PES General Meeting 2016.}
Frequency regulation service involves the injection or withdrawal of active power from the power grid to maintain the system frequency \cite{dvorkin_2014}. In the United States, frequency regulation is equivalent to secondary frequency control, while primary frequency control is more commonly known as frequency response. %~\cite{NERC_FC}.
The most common mean for a unit to provide frequency regulation is following a system operator's automatic generation control (AGC) signal, which computes the area control error (ACE) from frequency deviations and interchange power imbalances.

An appropriate response to the AGC signal requires dispatching flexibility within the assigned regulation capacity. Traditionally, the majority of frequency regulation capability has been provided by specially equipped generators. As technologies evolve, new types of flexibility resources emerge, such as battery and flywheel energy storage. These energy storage systems have significantly faster ramping capabilities compared to conventional generators, and their participation in frequency regulation can reduce the need to procure regulation capacity ~\cite{PNNL_Regulation}. To create an incentive for fast responding units to participate in frequency regulation, the Federal Energy Regulatory Commission (FERC) enacted FERC Order 755~\cite{FERC_755}, which requires system operators to add a performance payment with an accuracy adjustment to the capacity payment typically used in markets for ancillary services. This ``pay-for-performance" scheme has been implemented by most U.S. independent system operators (ISO) and regional transmission organizations (RTO). Under this new market rule, fast ramping regulation resources receive higher regulation payments because their performance is usually more accurate. A case study conducted by RES Americas has already shown that its 4~MW~/~2.6~MWh battery storage achieved a very high performance in following the regulation signal provided by PJM~\cite{RES_ES}.

To further utilize the fast responsive capability of energy storage systems (ESS) beyond the traditional AGC framework, some system operators, such as PJM, have introduced fast regulation. Units participating in fast regulation follow a regulation signal that changes much faster than the traditional AGC signal, and receive extra payments for doing so. The system operator also benefits from the improved regulation accuracy. A review by the California Energy Storage Alliance (CESA) shows that the combination of pay-for-performance and fast regulation has reduced the procurement requirement in the PJM regulation market by 30\%~\cite{CESA_Reg}.

However, the high investment cost of ESS has so far limited the installed energy storage capacity \cite{pandzic_2015}. The energy lost during charging and discharging and the energy imbalance in the regulation signal make it difficult for ESS operators to provide regulation services seamlessly over a relatively long time period (e.g. several hours). Some system operators have elected to lower this barrier by co-optimizing energy storage when dispatching regulation, i.e. taking into account the state of charge (SoC) and the cycling efficiency of ESS. 

This paper compares market policies that affect energy storage in the frequency regulation market. It covers the ISOs and RTOs that have organized wholesale electricity markets, including the PJM Interconnection (PJM); the New York Independent System Operator (NYISO); the Midcontinent Independent System Operator (MISO); ISO New England (ISO-NE); California Independent System Operator (CAISO). Study results on pay-for-performance, fast regulation, and storage energy compensation are discussed in Section~II, Section~III, and Section~IV, respectively. Section~V analyzes the regulation payment in different markets. Section~VI summarizes the findings.

\section{Pay-for-Performance}

All ISO~/~RTOs in our study have complied with FERC Order 755 and implemented  pay-for-performance in their regulation markets, except for NYISO and MISO which pay for regulation capacity only in their  day-ahead regulation markets (DAM).
A typical pay-for-performance regulation market involves the following terms:

\begin{itemize}
	\item \emph{Regulation capacity ($C$)} is the maximum and minimum regulation power that a unit offers during a given period. The regulation capacity is usually symmetrical so that a unit provides the same positive and negative capacity (i.e. $\pm 5$ MW). With the exception that CAISO separates regulation into up and down components. The capacity price in $\mathrm{\$/MW\cdot h}$, represents the option price to reserve 1~MW of capacity to provide regulation in the future. Not to be confused with the exercise price of the energy required to follow the regulation signal
	\item \emph{Regulation mileage ($M$)} is the sum of the absolute values of the regulation control signal movements. The mileage of a regulation resource output $\{P_1,\dotsc,P_n\}$ over a time period with $n$ signal steps and a regulation capacity $P_{max}$ is
	\begin{align}\label{Eq:mileage}
	M=\sum_{i=1}^{n}|P_i-P_{i-1}|/P_{max}%\quad.
	\end{align}
	 %$\sum_{i=1}^{n}|P_i-P_{i-1}|/P_{max}$. 
	 The unit for mileage is $\mathrm{\Delta MW/MW}$, representing the mileage of 1~MW of regulation capacity. The unit for mileage price is $\mathrm{\$/\Delta MW}$.
	\item \emph{Performance factor ($\rho$)} is a score from 0 to 1 that indicates a unit's performance in following the regulation signal. Most ISO/RTOs model the performance factor as a unit's accuracy in signal following (integration of absolute errors). PJM uses a more complex measure of performance that involves delay, correlation and accuracy~\cite{PJM_M12}.
\end{itemize}
Table~\ref{tab:pfp} shows the names used for these quantities in each market. 

\begin{table}[t]
	\begin{center}
		\centering
		% \footnotesize
		\caption{Pay-for-performance market terms}
		\label{tab:pfp}
		% \begin{footnotesize}
		\begin{tabular}{|c|c|c|c|c|}
			\cline{2-5}
			\multicolumn{1}{c|}{} & \multicolumn{2}{|c|}{Price}  & Regulation  & Performance  \Tstrut\\
			\multicolumn{1}{c|}{} &  Capacity &  Mileage & mileage &  factor \Bstrut\\
			\hline
			
			PJM 	& capability 	& performance 	& mileage 		& performance \Tstrut\\
			\cite{PJM_M11, PJM_M28}		& price 		& price			& 				& score\Bstrut\\
			\hline
			NYISO 	& capacity 		& movement 		& movement 	& performance\Tstrut\\
			\cite{NYISO_Bill}		& price		 	& price 		& 			& index\Bstrut\\
			\hline
			MISO 	& capacity 		& mileage 		& mileage 	& performance\Tstrut\\
			\cite{MISO_M5, chen2015development}		& price			& price		 	& 			& accuracy\Bstrut\\
			\hline
			ISO-NE  & capacity 		& service 		& movement 	& performance\Tstrut\\
			\cite{ISONE_Rule1}		& price			& price 		& 			& score\Bstrut\\
			\hline
			CAISO 	& capacity 		& mileage 		& mileage 	& accuracy\Tstrut\\
			\cite{CAISO_755}		& price			& price 		& 			& percentage\Bstrut\\
			\hline
			
		\end{tabular}
		%	\end{footnotesize}
	\end{center}
	\vspace{-.5cm}
\end{table}

In a pay-for-performance market, participants offer a capacity price bid and a mileage price bid. The system operator credits providers of regulation using the market capacity clearing price ($\lambda_C$), the mileage clearing price ($\lambda_M$), and the performance factor of its unit during the procurement period. A unit must reach a minimum performance score to be qualified for receiving credits, and maintain this minimum performance over a certain period to be eligible to bid into the regulation market in the future. We generalize the method for awarding regulation credit ($\lambda_R$) as follows: 
\begin{align}\label{Eq:reg_credit}
\lambda_R = C(\lambda_C+\rho M\lambda_M)%\quad.
\end{align}
Eq.~\ref{Eq:reg_credit} applies to CAISO and the real-time regulation market in NYISO. Eq.~\ref{Eq:reg_credit} also applies to MISO when a unit's performance is below 70\%, a unit receives full credit when above 70\% during each 5-minute performance test interval. At ISO-NE, the performance factor is used to adjust the capacity payment, thus Eq.~\ref{Eq:reg_credit} becomes % $C\rho(\lambda_C+M\lambda_M)$. 
\begin{align}\label{Eq:reg_creditA}
\lambda_R' = \rho C(\lambda_C+M\lambda_M)%\quad.
\end{align}
PJM uses mileage ratio in its credit calculations. The mileage ratio is the ratio between the mileage of the assigned regulation signal and the mileage of the traditional regulation signal ($M_{RegA}$, see Section~\ref{Sec:FR:PJM}). The formula for calculating the credit becomes
\begin{align}\label{Eq:reg_creditB}
\lambda_R'' = \rho C(\lambda_C+\frac{M}{M_{RegA}}\lambda_M)%\quad.
\end{align}
% In MISO, the implementation of pay-for-performance market is more conservative, because its performance measurement does not directly affect the calculation of the regulation credit. Instead, it is used as a threshold below which the payment that a resource receives is reduced.

\section{Fast Regulation}\label{Sec:FR}

A fast regulation signal is generated by applying a  high-pass filter to the AGC signal. Such signals have a fast ramping rate but are designed to have a zero-mean. They are therefore ideal for ESS interested in providing frequency regulation. Fast regulations have been implemented by PJM and ISO-NE, while MISO provides an alternative approach for utilizing fast responsive units.

A regulation signal's energy capacity requirement refers to the minimum energy capacity needed to follow perfectly the regulation signal, neglecting the efficiency and SoC constraints of an ESS. 
%The hourly energy requirement of a regulation signal is calculated by adding algebraically the energy charged and discharged from the storage at each time step over an hour. 
The energy requirement of the signal in that hour is the difference between the maximum and minimum SoC during the operating hour.
%The energy requirement of the signal in that hour is therefore the difference between the maximum and minimum energy levels. 
For a regulation signal $\mathbf{S}_h$ at hour $h$ with a step size of $T_s$ and $N$ steps per hour, the hourly energy requirement $E_{h}$ is therefore
\begin{align}\label{Eq:energy_offset}
E_{h} = \max_{j=1\dotsc N}\left(\sum_{i=1}^{j}T_sS_i\right)-\min_{j=1\dotsc N}\left(\sum_{i=1}^{j}T_sS_i\right)%\quad.
\end{align}

The energy balance of a regulation signal is the initial SoC level that allows a storage with the minimum energy capacity ($E_h$) to follow the signal perfectly, ignoring its efficiency and SoC constraints. A perfectly balanced signal (at 50\%) allows the SoC level of the energy storage to return to its initial level at the end of the operating period. Since Eq.~\ref{Eq:energy_offset} assumes that the initial SoC is zero (i.e. it simulates an ESS with negative capacity), the energy balance $\sigma_{h}$ for the signal $\mathbf{S}_h$ is
\begin{align}
\sigma_h = -\frac{1}{E_h}\min_{j=1\dotsc N}\left(\sum_{i=1}^{j}T_sS_i\right)%\quad.
\end{align}

\begin{figure}[t]%
\centering
\vspace{-.3cm}
\subfloat[Regulation energy hourly requirement (1MW capacity).]{
	\includegraphics[trim = 05mm 5mm 10mm 0mm, clip,  width = .9\columnwidth]{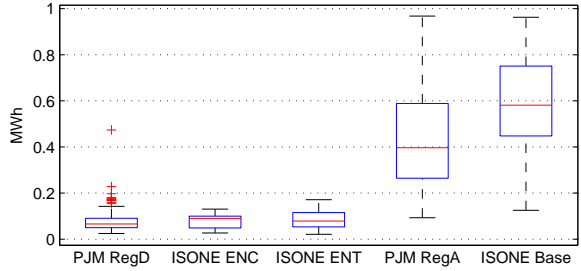}
	\label{Fig:energy_offset}%
}
% \\
\vspace{-.3cm}
\subfloat[Regulation energy hourly balance.]{
	\includegraphics[trim = 05mm 5mm 10mm 0mm, clip,  width = .9\columnwidth]{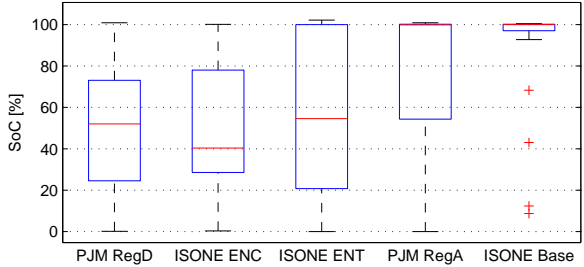}
	\label{Fig:energy_balancing}%
}
%\quad
\caption{Comparison of PJM and ISO-NE regulation signals.}%
\label{Fig:fast_reg}%
\vspace{-.4cm}
\end{figure}

\subsection{PJM RegD}\label{Sec:FR:PJM}

PJM provides two regulation signals in its day-ahead regulation market: RegA, is a low-pass filtered ACE designed for traditional regulation units, while RegD, is a high-pass filtered ACE designed for fast regulation units~\cite{PJM_M12}. The RegD signal is controlled to have zero-mean within 15 minutes~\cite{PJM_RegHist}. Units can choose freely which signal they follow, but can only follow one signal at a time. Our analysis of these regulation signals from June 2013 to May 2014 shows that the mileage of RegA is about 5~$\mathrm{\Delta MW/MW}$ per hour, while the mileage of RegD is around $15~\mathrm{\Delta MW/MW}$ per hour, suggesting that RegD ramps much faster. In market settlements, the mileage ratio for RegA is therefore 1, and the mileage ratio for RegD is around 3. Section~\ref{Sec:Reg_Pro} provides a detailed comparison of regulation payments. The box plot in Fig.~\ref{Fig:energy_offset} summarizes the distribution of the hourly energy requirement and shows that in most cases PJM RegD has an energy requirement of less than 0.25~MWh or approximately 15 minutes of full charge/discharge. Fig.~\ref{Fig:energy_balancing} shows that RegD has an average energy balance of 50\%, but that the extreme cases can reach 0 and 100\%.

\subsection{ISO-NE energy-neutral}\label{Sec:FF:ISONE}

ISO-NE includes ESS in the Alternative Technology Regulating Resources (ATRR), and provides two fast regulation signals (energy-neutral continuous~(ENC), energy-neutral trinary~(ENT)), and a conventional regulation signal~\cite{ISONE_EN_AGC}. A resource can change once a month the signal that it follows. ENC is identical to PJM's RegD signal, while ENT is specifically designed for ATRRs and the dispatch is either zero, full charge power, or full discharge power.
Since each trinary dispatched resource is assumed to have a step response and hence equal
participation factors, the ENT algorithm dispatches the entire set of Trinary Dispatched ATRRs together. 
Fig.~\ref{Fig:fast_reg} shows the energy requirement and the balance analysis using a one-day duration simulated signal from ISO-NE~\cite{ISONE_AGC} and suggest that these signals behave in a similar fashion as PJM's RegD and RegA . The mileage difference between the three simulated signals is not significant: the average mileages being 18, 16, and 14~$\mathrm{\Delta MW/MW}$ for ENC, ENT, and the conventional signal, respectively.

\subsection{MISO AGC-Enhancement}

In the MISO reserve market, the deployment interval for ESS (which are called stored energy resources or SER) is optimized for each ESS based on its initial condition~\cite{MISO_M2}. Taking things one step further, MISO is proposing a new regulation reserve deployment scheme called AGC-Enhancement~\cite{MISO_AGC}. The idea of this scheme is to deploy and un-deploy the fast ramping units first in the AGC dispatch, by creating a new fast ramping deployment group. MISO expects AGC-Enhancement to improve the utilization of fast ramping units.

\section{Storage Energy Offset Compensation}

Due to their less than perfect efficiency and the fact that the AGC signal has a non-zero mean, the SoC of an ESS is not maintained in the traditional regulation framework. While the fast frequency regulation alleviates  the effect of the signal energy imbalance, some ISOs have decided to directly compensate the energy offset of storage units in their regulation market through  dispatching. Such mechanisms include the Regulation Energy Management~(REM) of CAISO, and the Real-Time Dispatching~(RTD) of NYSIO. 

\subsection{CAISO REM~\cite{CAISO_REM}}

CAISO treats ESS as Non-Generating Resources~(NGR), a set of resources that also includes responsive demands. NGR units have two options for participating in the regulation market: non-REM~(traditional), and REM. Non-REM units are subjected to the same requirements as traditional regulation units, and must meet the 60-minute continuous energy requirement (i.e., to provide 1~MW regulation capacity, a non-REM NGR unit must be able to deliver at least 1~MWh). The continuous energy requirement for NGR-REM units is 15 minutes. NGR-REM units can only participate in the CAISO regulation market, and the subscription to REM can only be changed monthly. 

CAISO compensates the energy offset of NGR-REM units. By participating in REM, a NGR unit submits a preferred SoC set-point. CAISO will maintain this set-point by dispatching energy from the Real Time Market (RTM) for the next real time dispatch interval~\cite{CAISO_REMNGR}. This balancing energy is dispatched to each NGR-REM units together with the AGC signal by the Energy Management System (EMS). The EMS also takes into account a unit's efficiency in the dispatch. The REM energy compensation is guaranteed during normal operations, while in emergency cases the SoC is restored later.

In the CAISO regulation market settlement, the regulation energy charged and discharged from the storage and the dispatched REM energy are settled at the RTM locational marginal price (LMP), in addition to the capacity settlement and mileage settlement. Since REM guarantees to maintain the SoC of NGRs, the consequence of the energy settlement is that NGR-REM units must pay for their energy losses at their LMP. This settlement method is currently unique to CAISO.

% The CAISO regulation market separates regulation up and regulation down services, in regulation up, a unit only provide positive regulatory energy, and only negative in regulation down. A unit can provide both services at the same time, and the energy offsets in the two services are recorded separately. In the regulation market settlement, in addition to capacity and mileage settlements, the regulation up, regulation down, and the RTM offset energy are settled in locational marginal prices (LMP).

%A resource’s scheduling coordinator will allow the ISO to maintain the resource’s operating point by balancing the energy dispatched from the resource through the ISO Energy Management System to meet ISO regulation requirements. 

%The ISO will manage NGR-REM through the Real Time Market for the next Real Time Dispatch interval to offset the energy produced/consumed during the previous interval’s regulation energy dispatch.

\begin{figure}[t]
	\centering
	\vspace{-.5cm}
	\includegraphics[trim = 18mm 0mm 10mm 0mm, clip, width=.95\columnwidth]{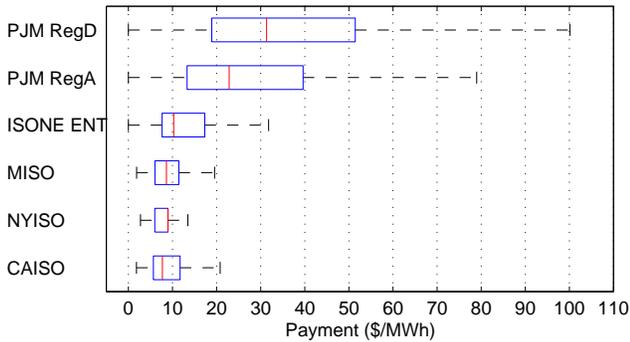}
	\caption{DAM regulation price comparison (outliers omitted).}
	\label{Fig:AS_price}
	\vspace{-5mm}
\end{figure} 

\subsection{NYISO RTD~\cite{NYISO_AS}, \cite{NYISO_LESR}}

NYISO calls ESS Limited Energy Storage Resources (LESR). LESR can bid into the regulation DAM and RTM, but only receive mileage payments in the RTM settlement. In the day-ahead regulation market, LESR units are evaluated and scheduled on an hourly basis without considering potential energy limitations. In the real-time regulation market, the RTD manages the energy level of LESR to maintain regulation capability to the extent possible by charging and discharging the LESR as necessary. The RTD assigns regulation base points to an LESR every  5 minutes based solely on its SoC level. For example, at each 5-minute interval, a 20~MW regulation capacity LESR  will provide $\pm$20~MW regulation reserve at 0~MW base point when its SoC is within a dead-band centered at 50\% SoC. When its SoC is above this dead-band, it will provide regulation with reduced capacity and a positive (charging) base point; when its SoC is below this dead-band, it will provide regulation with reduced capacity and a negative (charging) base point. Thus the SoC of the LESR converges to 50\%.

As in CAISO, the RTD SoC management is suspended during emergencies, when it is called RTD-CAM (Corrective Action Mode). During RTD-CAM, all LESR providing regulation discharge to zero in order to make up for the loss of generation. The SoC level is restored after RTD-CAM terminates.

% \section{Storage Energy Capacity Sizing }

\section{Regulation Market Income}\label{Sec:Reg_Pro}

\subsection{Day-ahead Market Results}

This section compares regulation payments in different regulation markets using public regulation pricing data from each ISO/RTO. In pay-for-performance markets, we assume that the performance factor is an ideal 100\%, and use Eqs.~\ref{Eq:reg_credit}-\ref{Eq:reg_creditB} to calculate these payments. The box plot in Fig.~\ref{Fig:AS_price} shows the distribution of calculated regulation payments in the DAM for each ISO/RTO.

\subsubsection{PJM}
The payments for the two regulation products available in PJM are based on signals and market clearing prices from June~2013 to May~2014. The hourly regulation mileage for each signal is calculated using Eq.~\ref{Eq:mileage}.

\subsubsection{CAISO}

CAISO publishes all its historical market data on the CAISO Open Access Same-Time Information System (OASIS). Besides pricing data, CAISO also publishes the hourly reg-up and reg-down mileage data (system mileage multipliers) of its AGC signal. To be consistent with the calculations for the other markets, we assume that each unit offers the same regulation up and regulation down capacity, we calculate the payments separately using day-ahead market prices from June~2014 to May~2015, and add the two results together. The average calculated reg-up price is \$5.2/MWh, while the reg-down price is \$4.1/MWh.

\subsubsection{MISO}

MISO publishes ancillary service pricing data on its market report website~\cite{MISO_data}. We use the day-ahead ancillary service market prices from June~2014 to May~2015. This market does not have a mileage price. 

\subsubsection{ISO-NE}

ISO-NE publishes ancillary service pricing data on its market data and information website~\cite{ISONE_data}. We use the final hourly clearing price data from June~2014 to May~2015. ISO-NE started mileage bidding on March 31st, 2015, so the clearing price for April and May contains a capacity price and a service (mileage) price. For these two months, we use mileage for the simulated ENT signal described in Section~\ref{Sec:FF:ISONE}, the payment for the other two signals are not calculated since the mileages are very close to those for ENT. 

\subsubsection{NYISO}

NYISO publishes its ancillary service pricing data on its pricing data website~\cite{NYISO_data}. We use the DAM regulation clearing prices from June~2014 to May~2015. This market does not have a mileage price.

% We use the regulation capacity and movement prices from June~2014 to May~2015 in its real-time ancillary service market, in which pay-for-performance is implemented. NYISO does not publish its AGC signal or mileage data, thus we use the AGC signal mileage data from CAISO to calculate the payments. 

%\subsection{Regulation Price Comparison between Markets}
%Fig.~\ref{Fig:AS_price} shows the payment comparison in different regulation DAMs. The regulation payments in PJM is significantly higher than the rest. The difference in regulation payment may due to extreme weather conditions, for example the regulation price in PJM of the year 2014 is about 45\% higher than 2013, primarily resulted from the record winter load~\cite{PJM_Report2014}. The structure of regulation markets, especially the ratio between the regulation capacity and the generation supply, also determines the price in regulation, this shall be discussed in detail in Section~\ref{Sec:Con}.

\subsection{Real-time Market}

CAISO, MISO, and NYISO have both RTM and DAM for operating reserves. The RTM is to procure additional reserves according to real-time system co-optimizations with updated system conditions~\cite{SANDIA_AS}. All three markets have mileage payments, however most mileage clearing prices are very close to zero and thus have little effect on the final price.
\subsubsection{CAISO}
CAISO determines its RTM regulation capacity in the hour-ahead scheduling process.
The weighted average RTM price for reg-up and reg-down is \$7.84/MWh and \$4.28/MWh, while for the RTM the prices are \$6.08/MWh and \$3.98/MWh, respectively~\cite{CAISO_Report2014}.
\subsubsection{MISO}
In MISO the RTM is cleared every five minutes to produce price and dispatch levels for regulation. Its annual average hourly RTM regulation prices are in general 10\% to 15\% higher than the DAM prices~\cite{MISO_Report2013}.

\subsubsection{NYISO}
NYISO calculates the RTM requirement over an RTD interval (1 hour) every 5~minutes. The average DAM regulation price is \$10.1/MWh in 2013 and \$12.9/MWh in 2014, while the average RTM regulation price is \$9.8/MWh in 2013 and \$13.8/MWh in 2014~\cite{NYISO_Report2014}. 

% The regulation payment in PJM significantly surpasses other markets, an immediate result of PJM's high regulation prices due to the success implementation of pay-for-performance and fast regulation. The payment in ISO-NE ranks the second, followed by MISO, NYISO, and CAISO at last. Notably, the mileage price in PJM is significantly higher than other markets, averaging around 4~$\mathrm{\$/\Delta MW}$, while the mileage prices in other market are very close to 0, mostly below 1~$\mathrm{\$/\Delta MW}$. 

\begin{table*}[t]
    \vspace{-12mm}
	\begin{center}
		\centering
		\footnotesize
		\caption{A Comparison of Regulation Market Policies on Energy Storage and Fast Responsive Units}
		\label{tab:rating}
		% \begin{footnotesize}
			\begin{tabular}{|c|c|c|c|c|c|c|c|c|c|}
\cline{2-9}
\multicolumn{1}{c|}{} & \multicolumn{2}{c|}{AS Requirements} &	ES Offset & Fast Reg. & ES & Reg. Cap. & Reg. Cost & Gen. Avg. \Tstrut\\
% \cline{2-4}
\multicolumn{1}{c|}{} & \multicolumn{2}{c|}{1~MW Cap.(MWh)} & Control & Schemes &  Abbr. & (MW) & (M\$) & (GW) \Bstrut\\
\hline
 %& 15 (REM) & & & & & & & \Tstrut\\
 %\cline{2-2}
CAISO \cite{CAISO_REM, CAISO_MO, CAISO_Report2014} & 0.25 (REM) & 1 (else) & REM & no &  NGR & $\sim325$ & 32  & $\sim28$ \Tstrut\Bstrut\\
\hline
PJM \cite{PJM_M11, PJM_M12, PJM_Report2014}& 0.25 (RegD) & 1 (RegA) & no & RegD &  - & $525-700$ & 254 & $\sim171$ \Bstrut\Tstrut\\
\hline
ISO-NE \cite{ISONE_Rule1, ISONE_OP14, ISONE_Report2014} & 0.25 (EN) & 1 (else)  & no & energy-neutral &  ATRR & $30-130$ &  29 & $\sim 12$ \Bstrut\Tstrut\\
\hline
MISO \cite{MISO_AGC, MISO_M2, MISO_Report2014} &  \multicolumn{2}{c}{-}  & \multicolumn{2}{|c|}{AGC-Enhancement}  & SER & $\sim 400$ & 20 & $\sim55$ \Bstrut\Tstrut\\
\hline
NYISO \cite{NYISO_AS, NYISO_Report2014} & \multicolumn{2}{|c|}{-} & RTD & no &  LESR & $150-250$ & 29 & $\sim 18$ \Bstrut\Tstrut\\
\hline

			\end{tabular}
	%	\end{footnotesize}
	\end{center}
	\vspace{-8mm}
\end{table*}

	\vspace{-1mm}
\section{Conclusions}\label{Sec:Con}
%	\vspace{-10mm}
	\vspace{-1mm}
% The participation of energy storages in U.S. regulation markets is becoming more friendly and more profitable, with the implementation of pay-for-performance markets, and the introduction of varies schemes on utilizing fast responsive units, and on compensating the energy offset of storages. 
Table~\ref{tab:rating} summarizes the findings of our study. In the column ancillary service (AS) requirements, we list the minimum energy capacity required to provide 1 MW of regulation reserve. This requirement is specified for CAISO, while for PJM and ISO-NE, it is derived from the simulations described in Section.~\ref{Sec:FR}. We were not able to find similar information for MISO and NYISO. The last three columns show the average required regulation capacity, the regulation cost, and the average generation capacity in each operator zone.

The comparison shows that PJM operates a very efficient regulation market in terms of the ratio between regulation capacity and generation, while the cost of regulation is similar to the other ISOs relative to the size of the interconnection. In particular, if we compare the ancillary market cost of PJM with CAISO, the historical cost of regulation in PJM per MWh load is $\sim$\$0.2-0.3/MWh, while this cost is $\sim$\$0.1-0.2/MWh for CAISO. On the other hand, the total per MWh ancillary service cost is much higher in PJM, at $\sim$\$2/MWh, than in CAISO, where it is $\sim$\$0.3/MWh~\cite{PJM_Report2014, CAISO_Report2014}. A possible explanation for this difference is that CAISO performs a co-optimization of energy and ancillary services, while PJM has a coupled co-optimization process with more complicated categories of ancillary services ~\cite{SANDIA_AS}. This latter approach is less efficient than CAISO's and results in a higher ancillary service cost. For example, sometimes the PJM RegD capacity is over-saturated, meaning that the capacity ratio of RegD is too high. Since the RegD logic is to maintain a 5-minute zero-mean energy, a high RegD capacity ratio can make the controller go against ACE correlation for extreme cases~\cite{PJM_RegHist}.

% From the comparison of regulation capacity and generation supply, we can clearly see that PJM has the most efficient regulation market, resulting in a much higher regulation price. The regulation price in PJM is about four times higher than CAISO, but the per MWh load cost is around \$0.2-0.3~/MWh, while this cost in CAISO is between \$0.1-0.2~/MWh. The ancillary service structure also contributes to the difference in regulation price. For example, in PJM the per MWh load cost of synchronized reserve is much lower than regulation, while in CAISO the cost of spinning reserve is higher than regulation, suggesting that the successful implementation of RegD has taken up the fast responsive task to frequency deviations~\cite{PJM_Report2014, CAISO_Report2014}.

Our comparison of different markets shows that the regulation markets are becoming more suitable and profitable for ESS because the new pay-for-performance policy makes energy storage more competitive. New zero-mean  regulation policies  have lowered the barrier for ESS to provide regulation for longer periods. From the perspective of a system operator, the participation of energy storage under appropriate market designs makes the regulation more efficient. However, providing a zero-mean regulation signal can be challenging, and must be well coordinated with other ancillary service products so that the energy balancing signal does not perturb the system.

\vspace{-2mm}
\section*{Acknowledgment}
\vspace{-1mm}
The authors thank Dr. I. Gyuk and his colleagues at the US DOE Energy Storage Program for funding this research. Sandia National Laboratories is a multi-program laboratory managed and operated by Sandia Corporation, a wholly owned subsidiary of Lockheed Martin Corporation, for the U.S. Department of Energy's National Nuclear Security Administration under Contract DE-AC04-94AL85000.
\vspace{-2mm}
\bibliographystyle{IEEEtran}	% (uses file "plain.bst")
\bibliography{IEEEabrv,literature}		% expects file "myrefs.bib"

\end{document}